\documentclass[11pt]{amsart2000}

\let\epsilon\varepsilon
\def\R{\mathbb{R}}
\def\N{\mathbb{N}}
\newtheorem{tw}{Theorem}
\newtheorem{cor}{Corollary}

\theoremstyle{definition}
\newtheorem{df}{Definition}

\newtheorem{ex}{Example}

\begin{document} 
\setlength{\unitlength}{0.01in}
\linethickness{0.01in}
\begin{center}
\begin{picture}(474,66)(0,0)
\multiput(0,66)(1,0){40}{\line(0,-1){24}}
\multiput(43,65)(1,-1){24}{\line(0,-1){40}}
\multiput(1,39)(1,-1){40}{\line(1,0){24}}
\multiput(70,2)(1,1){24}{\line(0,1){40}}
\multiput(72,0)(1,1){24}{\line(1,0){40}}
\multiput(97,66)(1,0){40}{\line(0,-1){40}}
\put(143,66){\makebox(0,0)[tl]{\footnotesize Proceedings of the Ninth Prague Topological Symposium}}
\put(143,50){\makebox(0,0)[tl]{\footnotesize Contributed papers from the symposium held in}}
\put(143,34){\makebox(0,0)[tl]{\footnotesize Prague, Czech Republic, August 19--25, 2001}}
\end{picture}
\end{center}
\vspace{0.25in}
\setcounter{page}{181}
\title[Lebesgue Theorem for multivalued functions]{On Lebesgue Theorem
for multivalued functions of two variables}
\author{Gra\.zyna Kwieci\'nska}
\address{University of Gda\'nsk\\
Institute of Mathematics\\
Wita Stwosza 57\\
80-952 Gda\'nsk, Poland}
\email{gkk@math.univ.gda.pl}
\thanks{This research was supported by the University of Gda\'nsk, grant
BW Nr 5100-5-0188-9}
\keywords{multivalued functions, semi-continuity of multivalued
functions, Baire classes of multivalued functions}
\subjclass[2000]{54C60, 54C08, 28B20}
\thanks{Gra\.zyna Kwieci\'nska,
{\em On Lebesgue Theorem for multivalued functions of two variables},
Proceedings of the Ninth Prague Topological Symposium, (Prague, 2001),
pp.~181--189, Topology Atlas, Toronto, 2002}
\begin{abstract}
In the paper we investigate Borel classes of multivalued functions of two
variables.
In particular we generalize a result of Marczewski and Ryll-Nardzewski
\cite{MN} concerning of real function whose ones of its sections are
right-continuous and other ones are of Borel class $\alpha$, into the case
of multivalued functions.
\end{abstract}
\maketitle

\section{Introduction}

Many results were publisched about the Borel classification of multivalued
functions depending on the one variable 
(see \cite{Kt, Ga, Br, Ha, Ma, Sl}). 
In the case of multivalued function of two variables we have the 
possibility of formulation of hypotheses concerning of its sectionwise
properties.

Lebesgue has shown that any real function $f$ of two variables with 
continuous ones of its sections and of Borel class $\alpha$ the other ones 
is of Borel class $\alpha +1$.
Marczewski and Ryll-Nardzewski have shown (see \cite{MN}) that the
condition of continuity in this theorem may be replaced by
right-continuity (or left-continuity).
In this paper we generalize these results into the case of multivalued
functions in possible general abstract spaces.
 
\section{Preliminaries}

Let $T$ and $Z$ be two nonempty sets and let $\Phi:T\to Z$ be a
multivalued function, i.e.\ $\Phi$ denotes a mapping such that $\Phi(t)$
is a nonempty subset of $Z$ for $t\in T$.
Then two inverse images of a subset $G\subset Z$ may be defined:
$$\Phi^+(G) = \{t\in T: \Phi(t)\subset G\}$$
and 
$$\Phi^-(G) = \{t\in T: \Phi(t)\cap G \ne \emptyset\}.$$

The following relations hold betwen these inverse images:
\begin{equation}
\Phi^-(G) = T\setminus \Phi^+(Z\setminus G)\ 
\mbox{and}\
\Phi^+(G) = T\setminus \Phi^-(Z\setminus G).
\end{equation}

Let $(T, \mathcal{T}(T))$ and $(Z, \mathcal{T}(Z))$ be topological spaces. 
The notations ${\rm Int}(A)$ and ${\rm Cl}(A)$ will be used to denote, 
respectively, the interior and the closure of a set $A$.

\begin{df} 
A multivalued function $\Phi:T\to Z$ is said to be $\mathcal{T}(T)$-upper
(resp.\ $\mathcal{T}(T)$-lower) semicontinuous at a point $t\in T$ if
$$
\forall G\in \mathcal{T}(Z) \ 
(\Phi(t)\subset G \Rightarrow t\in {\rm Int}\Phi^+(G))
$$
(resp.\ 
$\forall G\in \mathcal{T}(Z) \ 
(\Phi(t)\cap G\ne \emptyset \Rightarrow t\in {\rm Int}\Phi^-(G)))$.

$F$ is called $\mathcal{T}(T)$-continuous at the point t if it is 
simultaneously $\mathcal{T}(T)$-upper and $\mathcal{T}(T)$-lower
semicontinuous at $t$.

A multivalued function $\Phi$ being $\mathcal{T}(T)$-upper (resp.\
$\mathcal{T}(T)$-lower) semicontinuous at each point $t\in T$ is said to
be $\mathcal{T}(T)$-upper (resp.\ $\mathcal{T}(T)$-lower) semicontinuous.
\end{df}

It is clear that a multivalued function $\Phi$ is $\mathcal{T}(T)$-upper
(resp.\ $\mathcal{T}(T)$-lower) semicontinuous if and only if
$\Phi^+(G)\in \mathcal{T}(T)$ (resp.\ $\Phi^-(G)\in \mathcal{T}(T)$),
whenever $G\in \mathcal{T}(Z)$.

Given any countable ordinal number $\alpha$, let $\sum_{\alpha}(T)$ and
$\Pi_{\alpha}(T)$ denote the additive and multiplicative class $\alpha$,
respectively, in the Borel hierarchy of subsets of the topological space
$(T,\mathcal{T}(T))$.

We shall always assume $\alpha$ to be an arbitrary countable ordinal
number.

In perfect spaces the following inclusions hold:
\begin{equation}
\sum _{\alpha}(T)\subset \Pi _{\alpha+1}(T)\subset \sum _{\alpha+1}(T).
\end{equation}

\begin{df} 
A multivalued function $\Phi:T\to Z$ will be said to be of
$\mathcal{T}(T)$-lower (resp.\ $\mathcal{T}(T)$-upper) Borel class 
$\alpha$ if 
$$\Phi^-(G)\in \sum _{\alpha}(T)$$ 
(resp.\ $\Phi^+(G)\in \sum _{\alpha}(T)$), whenever $G\in \mathcal{T}(Z)$.
\end{df}

Let us note that a multivalued function of
$\mathcal{T}(T)$-lower (resp.\ $\mathcal{T}(T)$-upper) class $0$ is
$\mathcal{T}(T)$-lower (resp.\ $\mathcal{T}(T)$-upper) semicontinuous.

Let $f:T\to \R$ and $g:T\to \R$ be point-valued functions.
Then a multivalued function
$\Phi:T\to \R$ defined by formula
\begin{equation}
\Phi(t)=[f(t),g(t)]\subset \R
\end{equation}
is of $\mathcal{T}(T)$-lower (resp.\ $\mathcal{T}(T)$-upper) Borel class
$\alpha$ if and only if $f$ is of $\mathcal{T}(T)$-upper 
(resp.\ $\mathcal{T}(T)$-lower) and $g$ is of $\mathcal{T}(T)$-lower 
(resp.\ $\mathcal{T}(T)$-upper) class $\alpha$ in the Young
classification.

In fact, for $a<b$ we have
$$
\Phi^-((a,b))=\{t\in T: f(t)<b\}\cap \{t\in T: g(t)>a\}
$$
and
$$
\Phi^+((a,b))=\{t\in T: f(t)>a\}\cap \{t\in T: g(t)<b\}.
$$

\section{Main results}

Let $F:X\times Y\to Z$ be a multivalued function and
$(x_0,y_0)\in X\times Y$.
Then a multivalued function $F_{x_0}:Y\to Z$ such that 
$F_{x_0}(y)=F(x_0,y)$ is called $x_0$-section of $F$. 
Similarly a multivalued function $F^{y_0}:X\to Z$ such that
$F^{y_0}(x)=F(x,y_0)$ is called $y_0$-section of $F$.

\begin{tw}
Let $(Y,d)$ be a metric space and 
$(X,\mathcal{T}(X))$, $(Z,\mathcal{T}(Z))$ 
two perfectly normal topological spaces. 
Let $\mathcal{T}(Y)$ be a topology on $Y$ which is finer than the metric
one and such that $(Y,\mathcal{T}(Y))$ is separable. 
Let $S$ be a countable $\mathcal{T}(Y)$-dense subset of $Y$.
Suppose that to every point $v\in Y$ there corresponds a subset 
$U(v)\in \mathcal{T}(Y)$ such that
$$
\forall y\in S\ 
B(y)=\{v: y\in U(v)\}\in \sum _{\alpha}(Y,d)
$$
and
$$
\forall v\in Y\
\mathcal{N}(v)=\{U(v)\cap B(v,2^{-n}): n=1,2,\ldots\},
$$
where $B(v,2^{-n})$ denotes the open ball centered in $v$ with radius
$2^{-n}$, forms a filterbase of $\mathcal{T}(Y)$-neighbourghoods of the
point $v$.

Assume that $F:X\times Y\to Z$ is a multivalued function whose all
$y$-sections are of upper class $\alpha$ and all $x$-sections are
$\mathcal{T}(Y)$-continuous. 
Then $F$ is of lower class $\alpha +1$ on the product
$(X,\mathcal{T}(X))\otimes (Y,d)$.
\end{tw}

\begin{proof}
Let $D$ be an arbitrary $\mathcal{T}(Z)$-closed subset of $Z$. 
By (1) it is enough to show that 
$$F^+(D)\in \prod_{\alpha+1}((X,\mathcal{T}(X))\otimes (Y,d)).$$
Since $Z$ is perfectly normal, there is a sequence $\{G_n\}_{n\in \N}$ of
$\mathcal{T}(Z)$-open sets such that
\begin{equation}
D=\bigcap_{n \in \N}G_n=\bigcap_{n\in \N}{\rm Cl}(G_n)
\end{equation}
and
\begin{equation} 
{\rm Cl}(G_{n+1})\subset G_n \ 
\mbox{for}\ n\in \N.
\end{equation}
Let $S=\{y_k:k\in \N\}$. 
We will prove that
\begin{equation} 
F^+(D) = 
\bigcap_{n\in \N}\bigcup_{k\in \N}
(\{x:F(x,y_k)\subset G_n\}\times V_n(y_k)), 
\end{equation}
where
\begin{equation} 
V_n(y_k) = \{v\in Y:y_k\in U(v)\}\cap B(v,2^{-n}).
\end{equation}
Let 
$$(u,v)\in F^+(D)=\{(x,y)\in X\times Y:F(x,y)\subset D\}.$$
Then $F(u,v)\subset G_n$ for each $n\in \N$, by (4). 
Let $n$ be fixed. 
By the $\mathcal{T}(Y)$-upper semicontinuity of the $u$-section of $F$ at
the point $v\in Y$ there is a $\mathcal{T}(Y)$-open neighbourhood 
$U(v)\in \mathcal{N}(v)$ of $v$ such that $F(u,y)\subset G_n$ for any
$y\in U(v)$.

Let 
$$K=\{m\in \N:y_m\in U(v)\cap S\}$$ 
and let 
$$k={\rm min}\{m\in K:v\in V_n(y_m)\}.$$
Then 
$$(u,v)\in [F^{y_k}]^+(G_n)\times V_n(y_k)$$ 
and the inclusion
\begin{equation}
F^+(D) \subset 
\bigcap_{n\in \N}\bigcup_{k\in \N}
(\{x:F(x,y_k)\subset G_n\}\times V_n(y_k))
\end{equation}
is proved.

Conversely, let $(u,v)$ belongs to the right-hand side of (6). 
Suppose that $(u,v)\not \in F^+(D)$. Then by (4) we must have
\begin{equation}\label{q9}
F(u,v)\cap (Z\setminus {\rm Cl}(G_m)\neq \emptyset\ 
\mbox{for some}\ 
m\in \N.
\end{equation}
By $\mathcal{T}(Y)$-lower semicontinuity of the $u$-section of $F$ at the
point $v\in Y$ there is a $\mathcal{T}(Y)$-open neighbourhood 
$W(v)\in \mathcal{N}(v)$ of $v$ such that
\begin{equation}
F(u,y)\cap (Z\setminus {\rm Cl}(G_m)\neq \emptyset\
\mbox{for any}\ y\in W(v).
\end{equation}
We have supposed that 
$$
(u,v)\in \bigcap_{n\in \N}\bigcup_{k\in \N}
(\{x:F(x,y_k)\subset G_n\}\times V_n(y_k)).
$$
Therefore we conclude from (4) that to each $n$ there corresponds an 
index $k=k(n)$ such that
\begin{equation}
F(u,y_{k(n)})\subset G_n.
\end{equation}
For $v\in V_n(y_{k(n)})\subset B(v,2^{-n})$ we obtain 
$\lim_{n\to \infty}d(v,y_{k(n)})=0$. 
Since $y_{k(n)}$ tends to $v$ in $(Y,d)$ as $n$ tends to infinity,
(10) and (11) show that there is an index $n_0$ such that 
\begin{equation}
F(u,y_{k(n)})\cap (Z\setminus {\rm Cl}(G_m))\neq \emptyset\ 
\mbox{for any}\ n> n_0.
\end{equation}
By (5) and (11) we have 
$$F(u,y_{k(n)})\subset G_n\subset G_{n-1}\subset \ldots$$ 
for $n\in \N$.

In particular,
$$F(u,y_{k(n+j)})\subset G_{n+j}\subset G_n$$ 
for any $j\in \N$.
Fixing now $n=m$ (see (9)) we obtain $F(u,y_{k(m+j)}\subset G_m$ for any
$j\in \N$, which contradicts (12). 
We must have
$$
\exists {n\in \N} \ \forall {y\in S} \ 
v\not \in V_n(y)\vee F(u,y)\not \subset G_n.
$$
This formula means that
$$(u,v)\not \in \bigcap_{n\in \N}\bigcup_{k\in \N}
([F^{y_k}]^+ (G_n)\times V_n(y_k))$$
and the inclusion
\begin{equation}
\bigcap_{n\in \N}\bigcup_{k\in \N}
(\{x:F(x,y_k)\subset G_n\}\times V_n(y_k))\subset F^+(D)
\end{equation}
holds. 
By (8) and (13) the equality (6) is proved.

Observe that 
$$\{x:F(x,y_k)\subset G_n\}\in \sum_{\alpha}(X,\mathcal{T}(X))$$ 
since $y_k$-section of $F$ is of upper class $\alpha$. 
Furthermore it is assumed that $V_n(y_k)\in \sum_{\alpha}(Y,d)$. 
Therefore by (6) $F^+(D)$ is a countable intersection of countable unions
of the sets of the class 
$$\sum_{\alpha}(X,\mathcal{T}(X))\otimes \sum_{\alpha}(Y,d)\subset
\sum_{\alpha}(X\times Y),$$ 
where $X\times Y$ is the product of topological spaces 
$(X,\mathcal{T}(X))$ and $(Y,d)$. 
This completes the proof of Theorem 1.
\end{proof}

We give below two examples of topology $\mathcal{T}(Y)$ on $Y$ fulfilling
requirements of Theorem 1. 
From these examples it will be clear, that the $x$-sections of a 
multivalued function $F$ in Theorem 1 may be either all right-continuous
or all left-continuous in some meaning.

\begin{ex} 
Let $(Y,\diamond,d)$ be a topological group, whose topology is induced by
an invariant distance function $d$ (i.e.\ 
$d(\theta,y) = d(v, y\diamond v)$),
where $\theta$ denotes a neutral element of $Y$. 
Assume furthermore that $(Y,d)$ is separable.

Let $U\subset Y$ be an open set such that $\theta$ is an accumulation
point of $U$. 
Let 
$$U_n=(B(\theta,2^{-n})\cap U)\cup\{\theta\}
\mbox{ and }
V_n(y)=y\diamond U_n=\{y\diamond v:v\in U_n\}$$ 
for $n\in \N$. 
Then $\{V_n(y)\}_{n\in \N}$ forms a filterbase of neighbourhoods of a
point $y\in Y$ and the topology $\mathcal{T}(Y)$ in $Y$ generated by this
base fulfils all requirements of Theorem 1.

Indeed, it suffices to prove that $\{U_n\}_{ n\in \N}$ forms a base
of neighborhoods of $\theta$. 
We have
$$
U_n\cap U_m=U_{{\rm min}(n,m)}.
$$
Let $n\in \N$ and $v\in U_n$. 
Then there is $k\in \N$ such that
$$B(v,2^{-k})=v\diamond B(\theta,2^{-k})\subset U_n.$$
Therefore
$$
\forall n\in \N \ \forall v\in U_n \ \exists k\in \N \
V_k(v)\subset U_n.
$$

A countable dense subset of $(Y,d)$ is also $\mathcal{T}(Y)$-dense. 
It remains to show that $V_n(y)$ is a Borel set in $(Y,d)$ for any $n\in 
\N$.
Let $n\in \N$ and let $\Phi:Y\to Y$ be a multivalued function defined by
formula $\Phi(y)=V_n(y)$. 
Then $\Phi$ is continuous and and its graph
$$\operatorname{Gr}(\Phi)=\{(y,v):v\in \Phi(y)\}$$
is homeomorphic to the set 
$$Y\times U_n\in \sum_1(Y,d)\otimes (Y,d)\cap \prod_1(Y,d)\otimes (Y,d).$$
Finally $V_n(y)\in \sum_1(Y,d)\cap \prod_1(Y,d)$ for each $n\in \N$.
\end{ex}

\begin{ex} 
Let $(Y,d,\leq)$ be a linearly ordered metric space.
We follow Dravecky and Neubrunn (see \cite{DN}) in assuming that the 
space $(Y,d,\leq)$ has the property $\mathcal{U}$, i.e.\ $(Y,\leq)$ is
linearly ordered and there is a countable dense set $S$ in $(Y,d,)$ such
that for any $y\in Y$ we have $y=\lim_{n\to \infty}y_n$, where $y_n\in S$
and $y\leq y_n$ for $n\in \N$.
Then the topology $\mathcal{T}(Y)$ on $Y$ generated by all open sets in
$(Y,d)$ and also by all intervals $I_a=\{y\in Y:y\leq a\}$, $a\in Y$, 
fulfills the assumptions of Theorem 1. 
Indeed, let $y\in Y$ and $r>0$. 
Then
$$U_r(y)=B(y,r)\cap I_y=\{x\in Y:d(x,y)<r \wedge x\leq y\}$$
is a $\mathcal{T}(Y)$-neighbourhood of the point $y$.

Let $x\in U_r(y)$. 
Then $x\in B(y,r)$ and $x\leq y$, and then there is $r_1>0$ such that
$d(x,y)=r-r_1$. 
Let $\delta< {\rm min}(r,r_1)$. 
Then $B(x,\delta)\subset B(y,r)$. 
Let $n\in \N$ be such a number that $2^{-n}< \delta$. 
Then $U_{2^{-n}}(x)\subset U_r(y)$ and we see that 
$\{U_{2^{-n}}(y)\}_{n\in \N}$ forms a filterbase of
$\mathcal{T}(Y)$-neighbourhoods of the point $y$.

The set $S$ is also $\mathcal{T}(Y)$-dense. 
It remains to show that the set 
$$V_r(y)=\{z\in Y:y\in U_r(Z)\}$$ 
is a Borel set in $(Y,d)$. 
First we will show that
\begin{equation}
\begin{array}{l}
\mbox{If $y_0 \not = y$ and $y_0\in V_r(y)$},\ 
\mbox{then there exists}\ 0<r_1<r\\ 
\mbox{such that}\ 
U_{r_1}(y_0)\subset V_r(y)
\end{array}
\end{equation}
Suppose, contrary to our claim, that $U_{r_1}(y_0)\not \subset V_r(y)$ for
any $r_1<r$. 
Now let $n\in \N$ be such that $\frac{1}{n} < r$. 
Then there is $y_n$ such that $y\leq y_n$ and 
$y_n\in U_{\frac{1}{n}}(y_0)\setminus V_r(y)$, and then
$$
\begin{array}{lllllll}
y\leq y_n&
\wedge&
d(y_n,y_0)<{\frac{1}{n}}&
\wedge&
y_n \leq y_0&
\wedge&
(y_n\leq y \vee d(y_n,y)\geq r)
\end{array}
$$ 
for $n>{\frac{1}{n}}$.
If it were true that $d(y_n,y_0)<{\frac{1}{n}}$ and $y\leq y_n\leq y_0$ 
and $y_n\leq y$, we would have 
$$\lim_{n\to \infty}y_n=y_0=y,$$ 
in contradiction with $y \not = y_0$. 
Let $d(y_0,y)=\varepsilon$. 
If it were true that $d(y_n,y_0)<{\frac{1}{n}}$ and $d(y_n,y)\geq r$ we 
would have 
$$r\leq d(y_n,y)\leq d(y_n,y_0)+d(y_0,y)<{\frac{1}{n}}+\varepsilon.$$
Then we would have ${\frac{1}{n}}>r-\varepsilon>0$ for almost
every $n\in \N$, which is impossible. This establishes (14).

Our next claim is that
\begin{equation}
\begin{array}{l}
\mbox{If $y_0 \not =y$ and $y_0\in V_r(y)$}, 
\mbox{then there is}\ \delta>0\\ 
\mbox{such that}\ 
B(y_0,\delta)\subset V_r(y).
\end{array}
\end{equation}
Indeed, according to (14) there is $r_1\in (0,r)$ such that 
$U_{r_1}(y_0)\subset V_r(y)$. 
Let $\varepsilon=d(y_0,y)<r$ and let 
$\delta<{\rm min}(\varepsilon,r-\varepsilon,r_1)$.
Let $z\in B(y_0,\delta)$. 
Then either $d(y_0,z)<\delta$ and $z\leq y_0$ or $d(y_0,z)<\delta$ and
$y_0\leq z$. 
In the first case $z\in U_{\delta}(y_0)\subset V_r(y)$. 
In the second one 
$$d(z,y) \leq d(z,y_0)+d(y_0,y) < \delta+\varepsilon < 
r-\varepsilon+\varepsilon = r$$
and $y\leq z$ show that $z\in V_r(y)$. 
Combining these both results we conclude that 
$B(y_0,\delta)\subset V_r(y)$ and (15) is proved.

By (15) we see that the set 
$$\{z\in Y:d(z,y)<r \wedge y\leq z \wedge y\not =z\}$$
is open in $(Y,d)$. Therefore
$$
V_r(y) = 
\{y\}\cup \{z\in Y:d(z,y)<r \wedge y\leq z \wedge y\not =z\}\in 
\sum_1(Y,d)\cap \prod_1(Y,d).
$$

Note that this topology $\mathcal{T}(Y)$ may be viewed as a natural 
generalization of the known Sorgenfrey topology on the real line.
\end{ex}

\begin{cor} 
Let $f$ be a real function defined on the product of perfectly normal
toplogical space $X$ and the real line $\R$. 
Let us suppose that all $x$-sections of $f$ are right-continuous and all
$y$-sections of $f$ are of upper Young class $\alpha$. 
Then $f$ is of lower class $\alpha+1$ on $X\times \R$, i.e.\ it may be
represented as a point-limit of an increasing sequence of functions of 
upper Young class $\alpha$.
\end{cor}

\begin{proof}
Let us note that a multivalued function $F:X\times \R\to \R$ defined by 
formula
$$F(x,y)=[2-\arctan f(x,y),2+\arctan f(x,y)]$$
is of lower class $\alpha+1$, by Theorem 1. 
Moreover for $a<b$ we have
$$
F^-((a,b)) = 
\{(x,y): 2-\arctan f(x,y)< b\}\cap \{(x,y):2+\arctan f(x,y)> a\}.$$
By (3) the function $g(x,y)=2-\arctan f(x,y)$ is of upper class $\alpha+1$
and the function $h(x,y)=2+\arctan f(x,y)$ is of lower class $\alpha+1$ in
the Young classification, which finishes the proof of Corollary 1.
\end{proof}

The next theorem is a dualization of Theorem 1.

\begin{tw} 
Let $(Y,d)$ be a metric space and $(X,\mathcal{T}(X)),(Z,\mathcal{T}(Z))$ two
perfectly normal topological spaces. 
Let $\mathcal{T}(Y)$ be a topology on $Y$ which is finer than the metric one
and such that $(Y,\mathcal{T}(Y))$ is seperable. 
Let $S$ be a countable $\mathcal{T}(Y)$-dense subset of $Y$.
Suppose that to every point $v\in Y$ there corresponds a subset $U(v)\in
\mathcal{T}(Y)$ such that
$$
\forall y\in S\ B(y)=\{v: y\in U(v)\}\in \sum_{\alpha}(Y,d)
$$
and
$$
\forall v\in Y\
\mathcal{N}(v)=\{U(v)\cap B(v,2^{-n}): n=1,2,\ldots\},
$$
forms a filterbase of $\mathcal{T}(Y)$-neighbourghoods of the point $v$.
Let $F:X\times Y\to Z$ be a compact-valued multivalued function whose all
$y$-sections are of lower class $\alpha$ and all $x$-sections are 
$\mathcal{T}(Y)$-continuous. 
Then $F$ is of upper class $\alpha+1$ on the product 
$(X,\mathcal{T}(X))\otimes (Y,d)$.
\end{tw}

\begin{proof} 
Let $D$ be an arbitrary $\mathcal{T}(Z)$-closed subset of $Z$ and let
$S=\{y_k:k\in \N$. 
We will first prove that
\begin{equation}
F^-(D) = \bigcap_{n\in \N}\bigcup_{k\in \N}
(\{x:F(x,y_k)\cap G_n\ne \emptyset\}\times V_n(y_k)),
\end{equation}
where $G_n$ are open subsets of $Z$ fulfilling (4) and (5), while
$V_n(y_k)$ is defined by the formula (7).

If 
$$(u,v)\in F^-(D)=\{(x,y):F(x,y)\cap D\ne \emptyset\},$$ 
then by (4) $F(u,v)$ has nonempty intersection with $G_n$ for each $n\in \N$. 
Let $n$ be fixed and arbitrary. 
By $\mathcal{T}(Y)$-lower semicontinuity of $u$-section of $F$ at the 
point $v$ there exists a $\mathcal{T}(Y)$-open neighbourhood 
$U(v)\in \mathcal{N}(v)$ of $v$ such that $F(u,y)\cap G_n\ne \emptyset$
for all $y\in U(v)$.
Taking $k$ such that $v\in V_n(y_k)$ we have
$$(u,v)\in [F^{y_k}]^-(G_n)\times V_n(y_k) = 
\{x:F(x,y_k)\cap G_n\neq \emptyset\}\times V_n(y_k),
$$
which gives
$$
F^-(D) \subset 
\bigcap_{n\in \N}\bigcup_{k\in \N}
(\{x:F(x,y_k)\cap G_n\ne \emptyset\}\times V_n(y_k)).
$$
Now let us suppose that
$$
(u,v)\in \bigcap_{n\in \N}\bigcup_{k\in \N}
(\{x:F(x,y_k)\cap G_n\ne \emptyset\}\times V_n(y_k)).
$$
Then to each $n$ there corresponds an index $k=k(n)$ such that for
$y_{k(n)}\in S$ we have $F(u,y_{k(n)})\cap G_n\ne \emptyset$, and then by
(5)
\begin{equation}
F(u,y_{k(n+j)})\cap G_n\ne \emptyset\ 
\mbox{for any}\ j\in \N.
\end{equation}
If $(u,v)$ were not in $F^-(D)$, by (4) we would have
$$
F(u,v)\subset Z\setminus D=\bigcup_{n\in \N}(Z\setminus {\rm Cl}(G_n)).
$$
The value $F(u,v)$ is a compact subset of $Z$ and the sets $Z\setminus
{\rm Cl}(G_n)$, $n\in \N$, create a decreasing sequence of open sets, 
i.e.\ 
$$Z\setminus {\rm Cl}(G_n)\subset Z\setminus {\rm Cl}(G_{n+1}).$$
Therefore for some $m\in \N$ we have 
$F(u,v)\subset Z\setminus {\rm Cl}(G_m)$. 
Then by the $\mathcal{T}(Y)$-upper semicontinuity of $u$-section of $F$ at
the point $v\in Y$ we have $F(u,y)\subset Z\setminus {\rm Cl}(G_m)$ for
$y\in W(v)$, where $W(v)$ is a certain neighbourhood of the point $v$, 
chosen from the postulated filterbase $\mathcal{N}(v)$. 
Since $y_{k(n)}$ tends in $(Y,d)$ to $v$ as $n$ tends to infinity, by the
above there exists an index $n_0$ such that $y_{k(n)}\in W(v)$ for
$n>n_0$. 
Therefore
\begin{equation}
F(u,y_{k(n)})\subset Z\setminus {\rm Cl}(G_m)\
\mbox{for any}\ n>n_0.
\end{equation}
Taking $n=m$ in (17) we have 
$F(u,y_{k(m+j)})\cap G_m\ne \emptyset$ for any $j\in \N$,
which contradicts (18). Thus the equality (16) is proved.

Since the $y_k$-section of $F$ is of lower class $\alpha$, we have
$$\{x:F(x,y_k)\cap G_n\ne \emptyset\}\in \sum _{\alpha}(X).$$
Moreover under the assumption of our theorem we have 
$V_n(y_k)\in \sum_{\alpha}(Y,d)$. 
Thus we conclude from (16) that
$$
F^-(D)\in \sum_{\alpha}(X)\otimes \sum_{\alpha}(Y,d)\subset
\sum_{\alpha}(X\otimes Y)\subset
\prod_{\alpha+1}(X\otimes Y),
$$
where $X\otimes Y$ is the product of topological spaces
$(X,\mathcal{T}(X))$ and $(Y,d)$, as required.
The proof of Theorem 2 is finished.
\end{proof}


\begin{thebibliography}{1}

\bibitem{Br}
R.~Brisac, \emph{Les classes de {B}aire des fonctions multiformes}, C. R. Acad.
  Sci. Paris \textbf{224} (1947), 257--258. \MR{8,321f}

\bibitem{DN}
Jozef Draveck{\'y} and Tibor Neubrunn, \emph{Measurability of functions of two
  variables}, Mat. \v Casopis Sloven. Akad. Vied \textbf{23} (1973), 147--157.
  \MR{48 \#8735}

\bibitem{Ga}
K.~M. Garg, \emph{On the classification of set-valued functions}, Real Anal.
  Exch. (1985), no.~9, 86--93.

\bibitem{Ha}
Roger~W. Hansell, \emph{Hereditarily additive families in descriptive set
  theory and {B}orel measurable multimaps}, Trans. Amer. Math. Soc.
  \textbf{278} (1983), no.~2, 725--749. \MR{85b:54060}

\bibitem{Kt}
K.~Kuratowski, \emph{Some remarks on the relation of classical set-valued
  mappings to the {B}aire classification}, Colloq. Math. \textbf{42} (1979),
  273--277. \MR{81c:54024}

\bibitem{MN}
E.~Marczewski and C.~Ryll-Nardzewski, \emph{Sur la mesurabilit\'e des fonctions
  de plusieurs variables}, Ann. Soc. Polon. Math. \textbf{25} (1952), 145--154
  (1953). \MR{14,1070g}

\bibitem{Ma}
P.~Maritz, \emph{A note on semicontinuous set-valued functions}, Quaestiones
  Math. \textbf{4} (1980/81), no.~4, 325--330. \MR{83k:54016}

\bibitem{Sl}
W{\l}odzimierz {\'S}l{\c{e}}zak, \emph{Some contributions to the theory of
  {B}orel $\alpha$ selectors}, Problemy Mat. (1986), no.~5-6, 69--82.
  \MR{88h:54030}

\end{thebibliography}
\providecommand{\bysame}{\leavevmode\hbox to3em{\hrulefill}\thinspace}
\providecommand{\MR}{\relax\ifhmode\unskip\space\fi MR }
\providecommand{\MRhref}[2]{%
  \href{http://www.ams.org/mathscinet-getitem?mr=#1}{#2}
}
\providecommand{\href}[2]{#2}

\end{document}